\documentclass{article}
\usepackage{amsmath}
\usepackage{amssymb}
\usepackage{float}
\usepackage{graphicx}
\usepackage{natbib}
\usepackage{cite}
\usepackage[onehalfspacing]{setspace}
\usepackage{hyperref}
\usepackage{doi}
\title{Conditional Densities and Simulations of Inhomogeneous Poisson Point Processes: The R package "IPPP"}

\author{Niklas Hohmann\thanks{GeoZentrum Nordbayern, Fachgruppe Pal\"aoumwelt, University of Erlangen-Nuremberg, Loewenichstr. 28, 91054 Erlangen, Germany\; email:\href{mailto:niklas.hohmann@fau.de}{niklas.hohmann@fau.de}}}

\begin{document}

\maketitle

\begin{abstract}
A number of numeric approaches to simulate Poisson point processes with arbitrary event rates are presented and implemented for R. They include the simulation of the number of points and their location as well as the determination of conditional probability densities and random numbers in cases where additional information (e.g. location of points or number of points occurring) is available upfront.
\end{abstract}

\section{Introduction}
Inhomogeneous Poisson point processes (IPPPs) are versatile stochastic processes that describe points or "events" that occur with a changing rate and that are stochastically independent from each other. IPPPs can be used to model a variety of spatial and temporal  processes, such as the number of customers arriving at a supermarket as well as the time or the location of animal sightings.\\
The aim of this paper is to present a number of numerical approaches to generate random numbers corresponding to the location of points or the number of points in a set generated by an IPPPs. For the cases where additional information such as the locations of points or the number of points occurring is available upfront, procedures to derive the conditional probability densities are given in addition to the procedures to generate the corresponding random numbers.\\
IPPPs can be characterized in various ways. In the simple case, an IPPP on $\mathbb R^n$ is determined by a rate function $r\geq 0$, which determines rate of events \citep{daley2003}. In a more abstract setting, an IPPP on some localized space $S$ is characterized by an intensity measure $\mu$, assigning a bounded set $A$ the expected number of points $\mu(A)$ in $A$\citep{kallenberg2017}. The description via the rate function translates into the description via intensity measure by the relation
\begin{equation}
\mu(A)= \int_A r(x) \; \mathrm dx 
\end{equation}
Both characterizations are used in this manuscript to capture the full generality of the described methods.\\
For the case of IPPPs on $\mathbb R$ with rate functions, the described procedures are implemented for R\citep{R} and bundled in the package IPPP, which can be be downloaded on CRAN under \url{https://CRAN.R-project.org/package=IPPP }

\section{Inhomogeneous Poisson Point Processes on Bounded Sets}
Let $\zeta$ be an IPPP on $S$ with intensity measure $\mu$. Then on any bounded $A \subset S$ the following holds:\citep[theorem 3.4 and p. 72-73]{kallenberg2017}
\begin{itemize}
\item The number of points $\zeta(A)$ in $A$ follows a Poisson distribution with mean $\mu(A)$
\item On $ \{ \zeta(A)=l\}$, the locations $X_1, X_2, \dots, X_l$ of these points are i.i.d. with their distribution given by the probability measure $P( \; \cdot \; )= \mu( \; \cdot \;\cap A ) /\mu(A)$
\end{itemize}
Therefore the simulation if IPPPs on bounded sets decomposes into three steps: First the determination of $\mu(A)$, second the generation of a random number $l$ following a Poisson distribution with mean $\mu(A)$, and third the generation of $l$ random numbers with distribution $P$.\\
The case where the number of points occurring in $A$ is known to be a fixed number $m$ accordingly reduces to the generation of $m$ random numbers with distribution $P$.\\
In the case where $A \subset \mathbb R^n$ and the intensity measure $\mu$ is given by
\begin{equation}
\mu(B)=\int_B r(x) \; \mathrm dx  \quad  \text{ with } r \geq 0
\end{equation}
for a rate function $r$, rejection sampling with the restriction of the rate function $r$ onto the set $A$ is a convenient method to generate the random numbers following the distribution $P$. Since the rejection method does not require $r$ to be a probability density, no rescaling of $r$ is required. The determination of $\mu(A)$ by numeric integration can therefore be skipped in the case of conditioning on a fixed number $m$ of points occurring.\\
If additionally $n=1$, i.e. $\zeta$ is an IPPP on $\mathbb{R}$, and a fixed number $m$ of points occur at the (unordered) locations $X_1, X_2, \dots, X_m$, the points can be enumerated in increasing order and are then given by the order statistics $X_{(1)},X_{(2)}, \dots,X_{(m)}$. Let
\begin{equation}
f_P (x)= \begin{cases} r(x)  \left(\int_A r(x) \; \mathrm dx\right)^{-1} & \text{ for } x \in A \\ 0 & \text{ for } x \notin A \end{cases}
\end{equation}
be the probability density of $P$ and $F_P$ its cumulative distribution function. Then the probability density of the location of the $k$-th out of $m$ points is the probability density of the $k$-th order statistic and accordingly given by \citep[chapter 2]{david2004}
\begin{equation}
f_{k,m}(x)=k \binom mk [F_P(x)]^{k-1}[1-F_P(x)]^{m-k} f_P(x)\; ,
\end{equation}
More results related to this can be obtained from the theory of order statistics.

\section{Inhomogeneous Poisson Point Processes on the Real Numbers}
Let $\zeta$ be a homogeneous Poisson point process (HPPP) on $\mathbb R $ with rate 1, $\mu$ a locally finite measure on $\mathbb R$ and
\begin{equation}
R(t)= \begin{cases} \mu((0,t]) & \text{if } t\geq 0 \\ -\mu((t,0]) & \text{if } t < 0\end{cases}
\end{equation}
its distribution function. Then $\eta=\zeta \circ R$ is an IPPP with intensity measure $\mu$ \citep[theorem 24.16]{klenke2008}.\\
This allows to transform a realization of a HPPP into an IPPP with arbitrary intensity measure, just as inverse transform sampling transforms uniform distributed random numbers into random numbers with arbitrary distribution. The realization of a HPPP can be obtained by adding up exponentially distributed random numbers.\\
If a point $x_i$ of the IPPP is given, it can be transformed into $y_i=R(x_i)$. Treating $y_i$ as a point of the HPPP, the $n$-th point above/below $y_i$ can be generated by adding/subtracting a random number following an Erlang distribution with shape parameter $n$ and rate parameter 1 to/from $y_i$. Transforming this number back using $R^{-1}$ yields the $n$-th point above/below the known point $x_i$.\\
For the HPPP, the probability density $f$ of the $n$-th point above/below $y_i$ is given by the probability density of an Erlang distribution with shape parameter $n$ and rate parameter 1 (that is mirrored at the $y$-axis for points below $y_i$) and shifted by $y_i$. If $\mu$ is given by 
\begin{equation}
\mu(B)=\int_B r(x) \mathrm dx \quad  \text{ with } r \geq 0
\end{equation}
where $r$ is a rate function and $\lambda$ the Lebesgue measure, then the probability density $g$ of the $n$-th point of the IPPP above/below $x_i$ is given by 
\begin{equation}
g(x)=r(x) f(\pm R(x) -y_i)\; ,
\end{equation} 
see \citet[theorem 1.101]{klenke2008} in combination with \citet[theorem 5.2]{callahan2010}.

\section{Discussion}
With the upper procedures, the simulation od IPPPs in $\mathbb R$ with rate functions reduces to the standard problems of numeric integration, generating random numbers with some given distribution as well as the transformation of values using the distribution function or its inverse. For all of these problems, good solutions are available. The possibility to speed these solutions up strongly depends on the ability to exploit additional structures of the IPPP given. This can for example be the availability of closed representations of the distribution function or its inverse, or using a suitable procedure to generate the needed random numbers.

\end{document}